\documentclass[10pt,conference]{IEEEtran}
\usepackage[french]{babel}
\usepackage[T1]{fontenc}
\usepackage[latin1]{inputenc}
\usepackage{graphics, enumerate}
\usepackage{amsmath,amsfonts,amstext,amssymb,color,epsfig,multicol}

\graphicspath{{./figures/}}

\def\t{\theta}
\def\T{\Theta}
\def\s{\sigma}
\def\phi{\varphi}
\def\b{\beta}
\def\R{\mathbb{R}}

\def\IC{\text{IC}}

\def\ML{\text{ML}}
\def\AIC{\text{AIC}}
\def\BIC{\text{BIC}}
\def\IC{\text{IC}}

\def\1{ \hbox{ {\large 1} \!\!\!\!{\large I}} }
\def\qed{\hbox{$\vcenter{\vbox{
   \hrule height 0.4pt\hbox{\vrule width 0.4pt height 6pt
    \kern5pt\vrule width 0.4pt}\hrule height 0.4pt}}$}}

\begin{document}
\title{An improved method for model selection based on Information Criteria}

\author{
\authorblockN{Guilhem Coq}
\authorblockA{Laboratoire de Mathématiques et Applications (UMR 6086)\\
Université de Poitiers, Teleport 2, BP 30179 \\
86962 Chasseneuil Futuroscope Cedex\\
coq@math.univ-poitiers.fr}
\and
\authorblockN{Olivier Alata}
\authorblockA{Laboratoire Signal Image et Communication (FRE 2731)\\
Université de Poitiers, Teleport 2, BP 30179\\
86962 Chasseneuil Futuroscope Cedex\\
alata@sic.sp2mi.univ-poitiers.fr}
\and
\authorblockN{Marc Arnaudon}
\authorblockA{Laboratoire de Mathématiques et Applications (UMR 6086)\\
Université de Poitiers, Teleport 2, BP 30179\\
86962 Chasseneuil Futuroscope Cedex\\
arnaudon@math.univ-poitiers.fr}
\and
\authorblockN{Christian Olivier}
\authorblockA{Laboratoire Signal Image et Communication (FRE 2731)\\
Université de Poitiers, Teleport 2, BP 30179\\
86962 Chasseneuil Futuroscope Cedex\\
olivier@sic.sp2mi.univ-poitiers.fr}
}

\maketitle

\begin{abstract}
Information criteria are an appropriate and widely used tool for solving model selection problems. However, different ways to use them exist, each leading to a more or less precise approximation of the sought model. In this paper, we mainly present two methods of utilisation of information criteria : the classical one which is generally used and an alternative one, more precise but requiring a little more calculations. Those methods are compared on 1-D and 2-D autoregressive models ; we use a synthetized process for the 1-D case and texture images for the 2-D case. We also work with the original $\phi_\b$ criterion which includes all others usual criteria such as AIC, BIC, and $\phi$.
\end{abstract}

\section{Introduction}

An observation $x^n = x_1,\dots,x_n$ of a stochastic process $X$ and a parametric family of probability density functions $\{f(.|\t), \ \t \in \T\}$ being given, the Maximum Likelihood (ML) method allows to estimate a parameter $\hat \t \in \T$ fitting the observation. However, the problem of model selection is of greater interest. Let us cite for example the determination of the number of components of a mixture law, the order of an autoregression \cite{Hannan_Quinn,OO}, or of a Multiple Markov Chain \cite{Zhao_01}.

Unfortunately, for this problem, the ML method fails and overestimates the sought model. This is mainly due to the fact that there exists in $\T$ a parameter giving a high probability to the observation, even though that parameter may have many components. This is typically the case for an observation of length $n$ of a Multiple Markov Chain which may always be given a probability 1 if we suppose that its order is $n-1$.

An alternative method to ML is given by Information Criteria (IC). They are written under the general form $\IC = -\log(\ML) + \text{Pen}$, where $\text{Pen}$ is a penalty term growing as the parameter becomes complex. Since the term $-\log(\ML)$ has the opposite variation, the minimization of $\IC$ realizes a compromise between the data fitting and the complexity of the chosen parameter. Applications of those criteria are numerous, in signal processing as well as in pattern recognition \cite{OO}.

Different kinds of penalties are suggested. Based upon the minimization of a Kullback risk, Akaike \cite{Akaike_74} introduced the first criterion AIC ; Schwarz \cite{Schwarz} then suggested the BIC criterion using Bayesian estimation. Next, Rissanen used notions of coding and stochastic complexity \cite{Rissanen_86, Rissanen_89} to justify a criterion which has asymptotically the same expression as BIC. In the continuity of the work of Rissanen, El-Matouat and Hallin \cite{EMH} introduced the family of criteria $\phi_\b$. Note that the criterion $\phi$ given by Hannan and Quinn in \cite{Hannan_Quinn} is prior to $\phi_\b$ and is its limit case for $\b=0$. In a general frame, Nishii \cite{Nishii_88} gave sufficient conditions on the penalty for those criteria to be weakly or strongly consistent.  

In a first section, the problem of model selection is set, as well as the general method of utilisation of $\IC$ which requires too many computations. Subsections \ref{classique} and \ref{Nishii} describe the two methods we study: classical method and alternative method. The classical one, widely used, is based upon embedded models ; it has the advantage of requiring few computations but only gives a rough approximation of true model. The alternative one, referred to as ``Nishii method'', is presented by Nishii, Zhao and al. \cite{ZKB,Nishii_88,NKB} and allows a more significant selection of the model at the cost of slightly more computations. To our knowledge, this method is not often used but deserves attention. In section \ref{simulations} we compare the two methods in the case of 1-D or 2-D autoregressive models. Only the $\phi_\b$ criterion will be used since it includes AIC, BIC, and $\phi$ criteria.

\section{Model selection by IC}

Let $\{ \Omega^n ; {\cal A}^n ; f(.|\t),\ \t \in \T\}$ be a statistical structure, where $\T$ is a subset of $\R^m$ and  $x^n=x_1,\dots,x_n$ a realisation of the unknown density $f(.|\t)$. We choose a reference parameter $\t^0 = (\t^0_1,\dots,\t^0_m) \in \T$, usually the null vector. Let us denote by $S^\star$ the support of $\t$:
$$
S^\star = \{j \in [\![1,m]\!] \ | \ \t_j \neq \t^0_j\}
$$
where $ [\![1,m]\!] $ is the set of integers $\{1,\dots,m\}$. For any support $S$ we note $\T_S$ the set of parameters whose support is $S$.

Selecting the model is determining, from $x$, the support $S^\star$. Once a support $\hat S$ is chosen, the unknown parameter $\t$ is estimated in the ML sense in $\T_{\hat S}$.

Information Criteria are an appropriate tool for selecting the support. For $S \subset[\![1,m]\!]$, they have the general form:
\begin{equation}
\label{IC_S}
\IC(S) = - 2 \log f(x^n|\hat \t_S) + |S|\alpha(n)
\end{equation}
where $|S|$ is the cardinal of $S$ and $\hat \t_S$ is estimated in the ML sense in $\T_S$. The penalties $\alpha(n)$ for the criteria we use are:
\begin{equation}
 \begin{array}{ll}
  \label{pénalités}
  \text{-AIC criterion,}       &\alpha(n) = 2                \\
  \text{-BIC criterion,}       &\alpha(n) = \log n           \\ 
  \text{-$\phi_\b$ criterion,} &\alpha(n) = n^\b \log \log n\\
 \end{array}
\end{equation}
For a fixed $n$, adjusting the value of $\b$ in the penalty function (\ref{pénalités}) of the $\phi_\b$ criterion allows to obtain others criteria:
\begin{equation}\label{beta_AIC}
\begin{array}{rcl}
\b_{\AIC} &=& (\log 2 - \log\log\log n) / \log n       \\
\b_{\BIC} &=& (\log \log n - \log\log\log n) / \log n
\end{array}
\end{equation}
Consequently, we will only use the $\phi_\b$ criterion for $\b$ ranging from 0 to 1 ; $\b=0$ corresponds to the $\phi$ criterion. Moreover, in \cite{O_99} the following bounds on $\b$ are proposed:
\begin{equation}
\label{bornes_beta}
\b_{\min} = \frac{\log \log n}{\log n} \leq \b \leq 1-\b_{\min} = \b_{\max}
\end{equation}
It has been shown empirically in several contexts that, for a classic utilisation of IC (see section \ref{classique}), the value $\b_{\min}$ often gives the best results ; however the theoretical justification of this result has not been established by the authors yet. In our simulations, we present the value of $\b_{max}$ even though it gives poor results in most cases.

The selection of the support is then done via the minimization of $\IC(S)$ among all supports:
\begin{equation}
\label{méthode_bourrin}
\hat S = \text{Argmin} \{ \IC(S) \ | \ S \subset  [\![1,m]\!] \}
\end{equation}

A criterion is said strongly consistent if $\hat S$ converges
almost-surely (a.s.) to $S^\star$ as $n \rightarrow \infty$; it is
said weakly consistent if the convergence only is in
probability. Using the conditions of Nishii \cite{Nishii_88}, in
the case of a product statistical structure, the BIC and $\phi_\b$
criteria, $0<\b<1$, present a strong consistency. Those results are
extended to the linear regression model, including the autoregressive
models used here, by Nishii and al. in \cite{NKB}. Those conditions
hold with BIC and $\phi_\b$ criteria for the two methods we will
discuss: $\hat k$ defined by (\ref{IC_k} and \ref{méthode_classique}) converges a.s. to $k^\star$ and $\hat S$ defined by (\ref{méthode}) converges a.s. to $S^\star$.

The method (\ref{méthode_bourrin}) answers the problem of model selection, but requires many computations, see table I for details. Here, we study two lighter methods.

\section{The studied methods} \label{methodes}

\subsection{Classical method} \label{classique}

Let us take $m$ nested subsets of $\T$: $\T_1 \subset \dots \subset \T_m \subset \T$ called models of order $k\in [\![1,m]\!]$ ; for example $\T_k = \R^k$. 

The problem is then restricted to the determination, from $x$, of the order $k^\star$ of the smallest model $\T_{k^\star}$ containing the unknown parameter $\t$. To this end, we set
\begin{equation}
\label{IC_k}
\IC(k) = - 2 \log f(x^n|\hat \t_k) + |\T_k|\alpha(n)
\end{equation}
where $\hat \t_k$ is estimated in the ML sense in the model $\T_k$ and $|\T_k|$ is the number of free components of this model.

The selection of the order is done via the minimization of $\IC(k)$ among $k$:
\begin{equation}
\label{méthode_classique}
\hat k = \text{Argmin} \{ \IC(k) \ | \ k \in  [\![1,m]\!] \}
\end{equation} 

This method requires the least operations, see table I for details, but does not solve the problem of the determination of the support $S^\star$.

\subsection{Nishii method} \label{Nishii}

A reference parameter $\t^0 = (\t^0_1,\dots,\t^0_m) \in \T$ is fixed. Using the notation of (\ref{IC_S}), we set $\IC_{\text{ref}}=\IC([\![1,m]\!])$. This is the reference value of the criterion computed on the model $\T_m$ where all components are free. Then, for $j\in [\![1,m]\!]$, we set $\IC(-j) = \IC([\![1,m]\!] \backslash \{j\})$ the value of the criterion computed on the model where all components are free, except the $j$-th which is frozen to $\t^0_j$, generally 0. The Nishii method consists in choosing as an estimation of the support the set of indexes:
\begin{equation}
\label{méthode}
\hat S = \{ j   \in [\![1,m]\!] \ \vert \ \IC(-j) > \IC_{\text{ref}} \}
\end{equation}
Those are the important indexes in the sense that the criterion prefers the full model rather than the model where the $j$-th component is frozen. 

For a brief comparison of the different methods in terms of computations, let us suppose that each model of order $k \in [\![1,m]\!]$ in \ref{classique} has dimension $k$. The table I gives the number of operations required to solve the model selection problem, each computation of an IC being weighted by the dimension in which it has to be done, \emph{e.g.} 2 computations in dimension 5 count for 10 operations. $\\$

\begin{tabular}{|rccc|}
\multicolumn{4}{c}{Table I : comparison in terms of required operations} \\
\hline
Method:    & General (\ref{méthode_bourrin}) & Classical (\ref{méthode_classique}) & Nishii (\ref{méthode}) \\
Selection: & Support                         & Order                               & Support                \\
Operations:& $m2^{m-1}$                      & $m(m+1)/2$                          & $m^2$                  \\
\hline
\end{tabular}

$\\$

\section{Application in the autoregression case}\label{simulations}

Let us recall the expression of Gaussian autoregressive (AR) models in
$d$ dimensions:
\begin{equation}\label{eqar}
X_t = - \sum_{i \in S} a_i X_{t-i} + E_{S,t}
\end{equation}
where $S \in \mathbb{Z}_{d}$ is the set of indices associated to the regression, $E_S = (E_{S,t})_{t \in Z^d}$ is a Gaussian white noise with variance
$\s_S^2$.

\subsection{One-dimensional autoregression}\label{1DAR}

\subsubsection{Presentation}
In 1D, the classical used support $S$ of the model
is of the form $[\![1,k]\!]$ defining the model of order $k$, called
$\T_k$ (see \ref{classique}). As $\t_k =
\left\{ { {\bf a}_k, \s_k^2 } \right\}$ with ${\bf a}_k =
(a_1,\dots,a_k)$ and $\s_k^2$ is the variance of the associated
Gaussian white noise, $|\T_k| = k + 1$ while $|\T_S| = |S| + 1$. Selecting the order of the
model (see \ref{classique}) is finding $k$; while selecting the support for $\t^0 =
0$ (see \ref{Nishii}), is finding the indexes $j\in [\![1,m]\!]$ for which $a_j \neq 0$, $m$ being the maximum value of the order. 

The Yule-Walker equations allow to estimate the parameters in the ML sense and it is known that minus the maximal log-likelihood is equal to $n(\log (2 \pi \hat \s_S^2)+1)$. Dropping terms which do not depend on $k$ or $S$, the expression (\ref{IC_k}) and (\ref{IC_S}) of the criteria respectively become:
\begin{eqnarray*}
\IC(k) &=& n\log \hat \s_k^2 + k \, \alpha(n)   \\
\IC(S) &=& n\log \hat \s_S^2 + |S|\, \alpha(n) 
\end{eqnarray*}
where $\hat \s_k^2$ is the estimated variance assuming the model of
order $k$, and $\hat \s_S^2$ the one estimated supposing the support
is $S$. A realisation $x^n$ of that process being given, we may apply
the two methods (\ref{méthode_classique}) and (\ref{méthode})
discussed above. Typically, if $ {\bf a} =(-1,0,1)$, we expect the classical method to choose order $\hat k =3$ and the Nishii method to choose support $\hat S= \{1,3\}$.

We generate 100 observations $x^n$ of an AR process (\ref{eqar}) of order 15 and parameters 
$$
{\bf a} = (0.5,\,0.4,\,0,\,0,\,0,\,0,\,0,\,0,\,0,\,0,\,0,\,0,\,0,\,0,\,0.45), \ \s^2 = 1
$$
and for each of these observations, we solve the model selection
problem using both classical and Nishii method with the $\phi_\b$
criterion. We set our maximal order $m$ to $20$. The classical method
is a success if it chooses $\hat k = 15$, while the Nishii method is a success if it chooses $\hat S = \{1,2,15\}$. 

\subsubsection{Results and discussion}

\begin{figure}[h]
\centering
\includegraphics[width=3in]{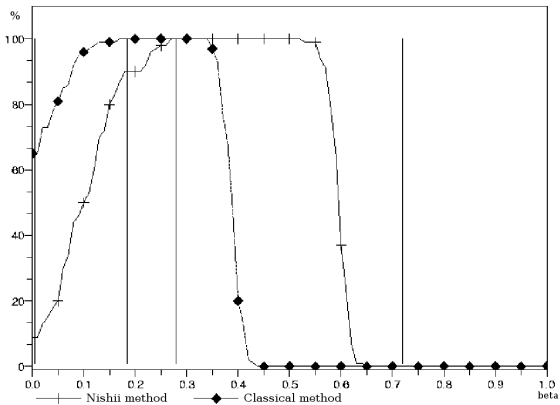}
\caption{Percentage of success for both methods, n=1000}
\label{fig:pourc_1000}
\end{figure}

Figure \ref{fig:pourc_1000} shows the percentage of succes of each method for $n=1000$. The $x$ axis represents the value of $\b$ used in the $\phi_\b$ criterion. The vertical lines correspond to the value of $\b_{\AIC}, \ \b_{\BIC}, \b_{\min}$ and $\b_{\max}$, always in that order ; see equations (\ref{beta_AIC}) and (\ref{bornes_beta}).

\begin{figure}[h]
\centering
\includegraphics[width=3in]{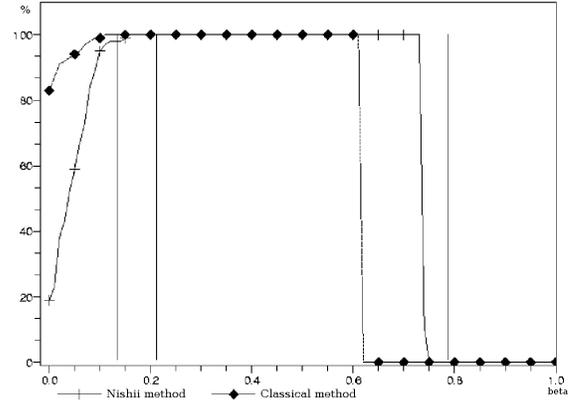}
\caption{Percentage of success for both methods, n=100000}
\label{fig:pourc_100000}
\end{figure}

We note that the AIC criterion often fails, especially with the Nishii method. The BIC criterion sometimes fails with the Nishii method, but the $\phi_{\b_\text{min}}$ criterion gives $100\%$ success with both methods. For small values of the penalty \emph{i.e.} $\b$ close to 0, IC gets close to the ML method, thus overparametrize the model. Moreover the Nishii method is less efficient in this area because if it keeps just one index in $[\![3,14]\!] \cup [\![16,20]\!]$, it fails ; while the classical method only fails if it chooses an order $\geq 16$. By opposition, for strong values of the penalty, IC tend to underparametrize the model. This happens here for the classical method and $\b \approx 0.45$: it only chooses order 2, thus misses $a_{15}=0.45$. The same happens for the Nishii method, but for $\b \approx 0.65$: it chooses support $\hat S = \{1,15\}$, thus misses the parameter $a_2=0.4$ which is the smallest. For $\b$ close to 1, both methods underparametrize so much that they choose to keep no parameter at all. The same results are presented for $n=100\ 000$ in figure \ref{fig:pourc_100000},  note that  $\b_{\AIC} <0$ as soon as $n \geq 1619$.

\begin{figure}[h]
\centering
\includegraphics[width=3in]{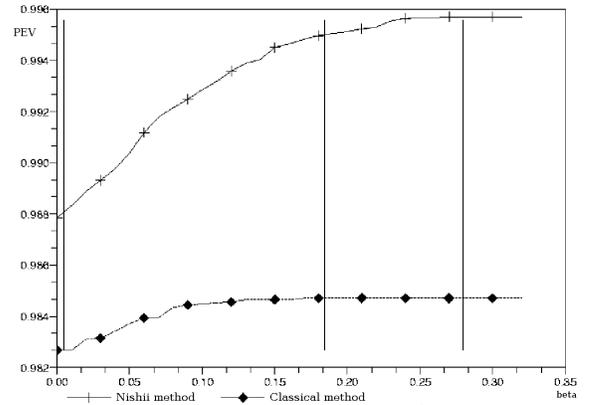}
\caption{Prediction error variance, n=1000}
\label{fig:erreur}
\end{figure}

Figure \ref{fig:erreur} presents the prediction error variance (PEV) of the models chosen by both methods for $0 \leq \b \leq 0.35$, \emph{i.e.} before the classical method starts to underparametrize. 

The more parameters are kept, the better the model fits the data, the smaller is the PEV. This explains why the PEV grows with $\b$ and why it is greater with Nishii method in the 100$\%$ success zone: Nishii method sets $a_3=\dots=a_{14}=0$ while the classical method estimates them. However, PEV with the Nishii method is closer to the real one $\s^2 = 1$. In that sense, the Nishii method appears to describe the model more precisely and the minimization of the PEV, equivalent to the ML method here, should not be a guideline for model selection.

\begin{figure}[h]
\centering
\includegraphics[width=3in]{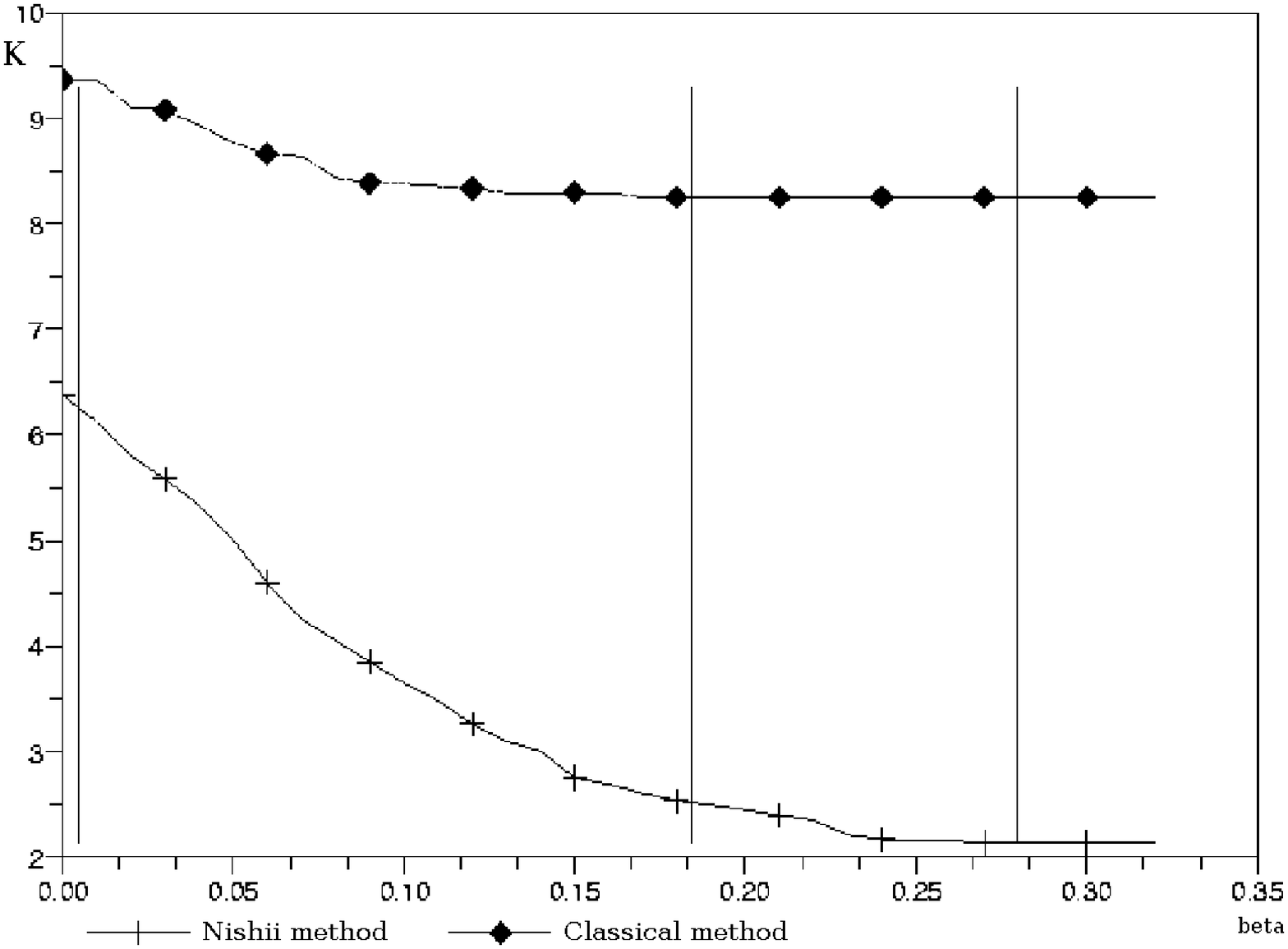}
\caption{Kullback distance to the true model, n=1000}
\label{fig:Kullback}
\end{figure}

Figure \ref{fig:Kullback} shows for the same values of $\b$ the Kullback distance between the true model $(a,\s)$ and the chosen one~$(\hat a, \hat \s)$:
$$
K\left((a,\s);(\hat a, \hat \s)\right) = - \frac{n}{2} + \log \frac{\hat \s}{\s} + \frac{\s^2}{2\hat \s^2}\text{Tr}\left({}^t(\hat A A^{-1})(\hat A A^{-1})\right)
$$
where $A$ and $\hat A$ are $n \times n$ matrix depending on $a$ and $\hat a$ respectively:
$$A=
\left(
\begin{array}{ccccccc}
    1  & 0      & \cdots & \cdots & \cdots & \cdots & 0      \\
  a_1  & \ddots & \ddots &        &        &        & \vdots \\
\vdots & \ddots & \ddots & \ddots &        &        & \vdots \\ 
  a_k  &        & \ddots & \ddots & \ddots &        & \vdots \\
   0   & \ddots &        & \ddots & \ddots & \ddots & \vdots \\
\vdots & \ddots & \ddots &        & \ddots & \ddots & 0      \\
   0   & \cdots &    0   &  a_k   & \cdots &    a_1 & 1      \\
\end{array}
\right)
$$

The Nishii method is seen to give a better description of the sought model in terms of Kullback distance.

\subsection{Two-dimensional autoregression}
\subsubsection{Presentation}
The support of the 2D AR model now contains couples of
integers. In litterature, the classical approach is based on
supports of different types of geometry \cite{Ranganath_85}: causal Quarter Plane
(QP), causal Non-Symetrical Half Plane (NSHP), semi-causal or
Non-Causal (NC). As 2D spectrum estimation methods based on QP support
provide nice results \cite{Alata_97}, we used here this type of
support.

Around a site, four QP supports can be defined. But, due to central symmetry, only two QP are
associated with different sets of AR parameters. The first one is, with order $(k_1,k_2)$:

\begin{equation*}\label{QP1}
QP1_{k_1,k_2} = \left\{  (i_1,i_2) \left| \begin{array}{c}
                                      0 \leq i_1 \leq k_1, 0 \leq i_2 \leq k_2 \\
                                      (i_1,i_2) \neq (0,0) 
                                     \end{array}
                                   \right.
                \right \}
\end{equation*}
while the second QP is :
\begin{equation*}\label{QP2}
QP2_{k_1,k_2} = \left\{ (i_1,i_2) \left| \begin{array}{c}
                                          -k_1 \leq i_1 \leq 0, 0 \leq i_2 \leq k_2 \\
                                          (i_1,i_2) \neq (0,0)  
                                         \end{array}
                                  \right.
                \right\}
\end{equation*}

The classical 2D QP AR model of order $(k_1,k_2)$ is:
\begin{equation*}\label{eqar2D}
X_{t_1,t_2} = - \sum_{(i_1,i_2) \in QP_{k_1,k_2}} a_{i_1,i_2} X_{t_1-i_1,t_2-i_2} + E_{QP,t_1,t_2}
\end{equation*}
where QP is either QP1 or QP2. We define $\T_{k_1,k_2}$ as the set of
parameters of 2D QP AR model of order $(k_1,k_2)$ so that $|\T_{k_1,k_2}| =
(k_1 + 1) \times (k_2 + 1)$, adding the variance of the prediction
error to the set of AR parameters.

By opposition to the
Nishii method which works as in the 1D case (each parameter associated
with a couple of integers
can be tested equal or not to zero), the increment in the
cardinality of nested models is not always one. For example,
$\T_{k_1,k_2+1}$ and $\T_{k_1+1,k_2}$ contains respectively $(k_1+1)$
and $(k_2+1)$ more parameters than $\T_{k_1,k_2}$. This fact implies that some indexes can be rejected by the classical method even if one of them would have been kept by the Nishii method.   

\subsubsection{Results and discussion} 
\begin{figure}[h]
\begin{center}
\begin{tabular}{c}
\includegraphics[width=0.25\textwidth]{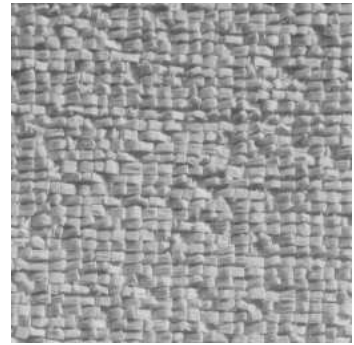} \\
(a) d84 texture\\
\includegraphics[width=0.25\textwidth]{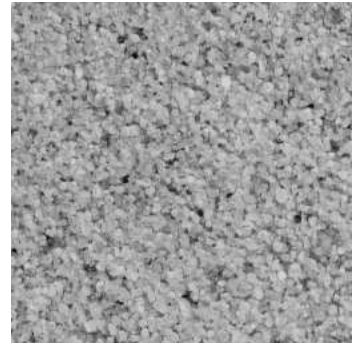} \\
(b) d29 texture \\
\end{tabular}
\end{center}
\caption{256$\times$256 textures from the Brodatz album}
\label{fig:textures}
\end{figure}

For running simulations, we used two textures from the Brodatz album
\cite{Brodatz_66} (see Figure \ref{fig:textures}) in order to show the
application of the Nishii method on real 2D processes. 

\begin{figure}[h]
\centering
\includegraphics[width=3in]{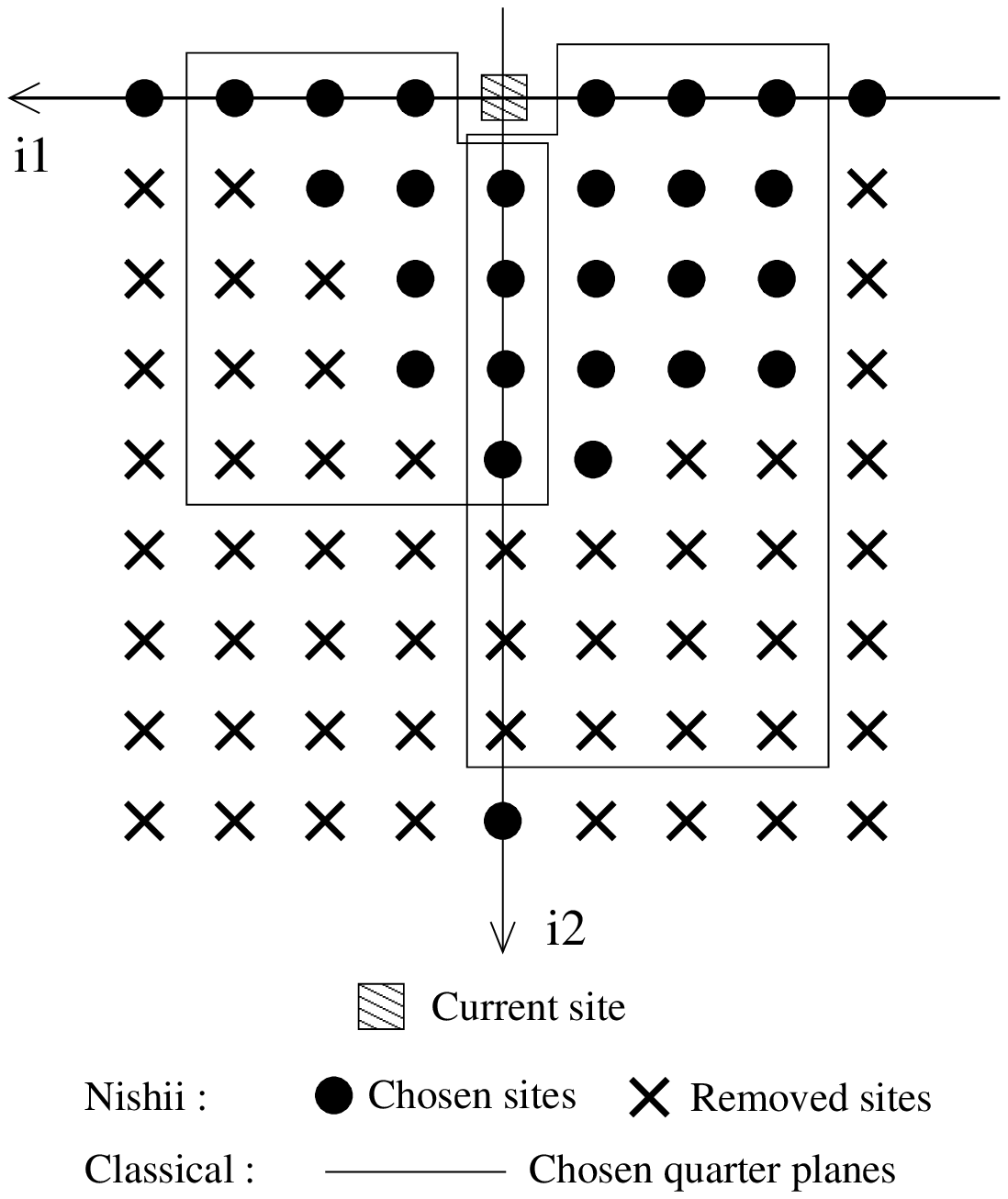} \\
\caption{Results of classical and Nishii methods on d84 texture}
\label{fig:qp_d84}
\end{figure}

\begin{figure}[h]
\centering
\includegraphics[width=3in]{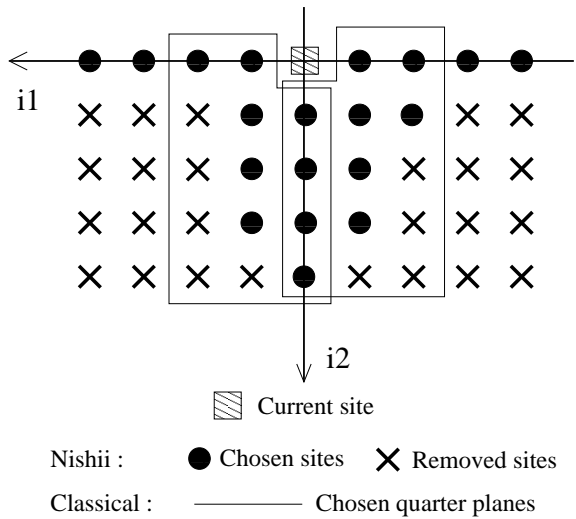} \\
\caption{Results of classical and Nishii methods on d29 texture}
\label{fig:qp_d29}
\end{figure}

We set our maximal order to $(m_1,m_2) = (18,18)$ and use classical and Nishii methods together with $\phi_{\b_{\min}}$ criterion for determining respectively the order and the support of the autoregression. Figures \ref{fig:qp_d84} and \ref{fig:qp_d29} present the results, on the left of the current site is QP1, on the right is QP2. 


Since it has to select rectangular supports, the classical method keeps sites which are not considered important by the Nishii method. Conversely, as noted earlier, the Nishii method keeps sites which are missed by the classical one. In the 1D synthetized case, we saw in figure \ref{fig:Kullback} that the Nishii method gives a more precise description of the model. Here, even though we did not suppose that our observation effectively comes from a true model, the model selected by the Nishii method is still more accurate. Moreover, as a perspective, the shape of the supports chosen by the Nishii method might be a discriminating factor between different texture images which might be used, for example, to improve recognition methods.


%

\end{document}